\documentclass[11pt,a4paper]{article}
\usepackage[T1]{fontenc}
\usepackage[utf8]{inputenc}
\usepackage{lmodern,microtype}
\usepackage{amsthm}
\usepackage{graphicx}
 \usepackage{booktabs}
\usepackage{cancel}
\usepackage{graphicx,subcaption}
\theoremstyle{definition}
\usepackage{comment} 
\usepackage{ulem}
\newtheorem{example}{Example}
\usepackage[margin=27mm]{geometry}
\usepackage{amsmath,amssymb,mathtools}
\usepackage[round,authoryear]{natbib}
\usepackage{xcolor}
\usepackage{placeins}
\usepackage[colorlinks=true,linkcolor=blue!40!black,citecolor=blue!40!black,urlcolor=blue!40!black]{hyperref}
\newtheorem{theorem}{Theorem}[section]

\newtheorem{corollary}[theorem]{Corollary}
\theoremstyle{definition}
\newtheorem{assumption}[theorem]{Assumption}
\newtheorem{definition}[theorem]{Definition}
\usepackage{bookmark}
\title{Anytime-valid inference for time series via self-normalization}
\author{Tim Kutta \and Xiaofeng Shao}

\allowdisplaybreaks[2]

\begin{document}

\maketitle

\begin{abstract}
\noindent
We develop a self-normalized approach to constructing asymptotic
confidence sequences for a broad class of parameters of stationary
time series under weak dependence. Our theory covers nonlinear
marginal parameters, including quantiles, expectiles and expected
shortfall, as well as autoregressive coefficients. To the best of our knowledge, this is the first use of self-normalization to construct asymptotic confidence sequences for time series.\\
A key advantage
is that the construction avoids direct estimation of the long-run
variance, which can be particularly challenging for nonlinear
parameters under serial dependence. Following an initial burn-in
period of length $m$, confidence intervals are released sequentially,
with simultaneous coverage over all subsequent times converging to
$1-\alpha$ as $m\to\infty$. We propose two sequential self-normalizers
with complementary statistical and computational advantages. The
forward construction uses only the history of sequential parameter
estimates, leading to minimal computational overhead. The three-block construction typically yields narrower confidence sequences
and is compatible with the optimal law-of-the-iterated-logarithm width rate. Together,
these constructions extend the benefits of self-normalization to
sequential inference under dependence. We investigate their finite-sample performance through simulations
and illustrate their use in Markov chain Monte Carlo output analysis
in the supplementary appendix.
\end{abstract}
\section{Introduction}\label{sec:introduction}

In time-series analysis, data are often collected sequentially, with
parameter estimates and confidence intervals updated with every new
observation. The aim of such an analysis is to accumulate evidence in
real time and draw conclusions when a suitable stopping criterion
is met. Stopping criteria may differ across applications.
Examples might be excluding a certain parameter value from a sequential
confidence interval (and rejecting the associated hypothesis),
reaching a prescribed interval width, or simply the analyst deciding
to end the experiment based on the observations collected so far.
Conventional confidence intervals are designed for a fixed,
predetermined sample size and generally do not preserve their coverage
guarantees under repeated peeking or data-dependent stopping.
In contrast, confidence sequences provide simultaneous coverage over
time, allowing the analyst to look at confidence intervals at any time
and to use data-dependent stopping rules
\citep{darling:robbins:1967,howard:ramdas:mcauliffe:sekhon:2021}.\\
Coverage control is particularly difficult for time series, where
interval widths must account for unknown serial dependence. This
difficulty is familiar from retrospective analyses, where conventional
confidence intervals often rely on long-run variance estimates. For
nonlinear parameters, such as quantiles, obtaining these estimates can
be especially tricky because it can require estimating some features
of the marginal distribution, such as its density at the quantile of
interest \citep{shao:2010}. A well-established approach to avoiding
direct long-run variance estimation is self-normalization. Roughly
speaking, this approach uses fluctuations of subsample estimates to
construct a suitable normalization \citep{lobato:2001,shao:2015}.
We will explain self-normalization in the retrospective setting in
greater detail below, as this provides the basis for understanding the
different ways in which it can be generalized to sequential inference.\\
To formulate our inferential objective more mathematically, let
$(X_t)_{t\in\mathbb N}$ be a stationary time series with marginal
distribution $F$, and let $\theta\in\mathbb R$ denote a parameter of
interest. For marginal parameters, we write $\theta=\theta(F)$.
For an initial burn-in period of length $m$, let $C_{n,m}$ be a
confidence interval based on $X_1,\ldots,X_n$, available at each time
$n\ge m$. We call $(C_{n,m})_{n\ge m}$ an asymptotic
$(1-\alpha)$-confidence sequence if
\begin{equation}\label{eq:intro-asymptotic-cs}
\liminf_{m\to\infty}
\mathbb P\!\left(\theta\in C_{n,m}\ \text{for all }n\ge m\right)
\ge 1-\alpha,
\qquad \alpha\in(0,1).
\end{equation}
This is the notion of asymptotic time-uniform coverage considered by
\citet[Definition~2.7]{waudby-smith:arbour:sinha:kennedy:ramdas:2024}.
We are interested in asymptotically sharp versions, where the coverage
probability converges to exactly $1-\alpha$.\\
For i.i.d. data, various papers have developed nonasymptotic confidence
sequences, with coverage guarantees that hold simultaneously for all
sample sizes. \citet{howard:ramdas:mcauliffe:sekhon:2021} use
concentration inequalities and nonnegative supermartingales to
construct confidence sequences for means and related quantities, with
applications to average treatment effects and covariance matrices.
\citet{howard:ramdas:2022} develop confidence sequences for quantiles
of i.i.d. observations, including guarantees that hold simultaneously
over time and over all quantile levels. Asymptotic constructions have
been considered, for example, by
\citet{waudby-smith:arbour:sinha:kennedy:ramdas:2024}, who use strong
approximation arguments to obtain confidence sequences for means,
with extensions to asymptotically linear estimators under i.i.d.
sampling and applications to average treatment effects.
\citet{bibaut:kallus:lindon:2024} establish asymptotic error control
and expected-rejection-time guarantees for sequential tests and obtain
associated confidence sequences by test inversion. Their framework
accommodates parameters defined through estimating equations, with the
mean and average treatment effect as examples. Another relevant point
of comparison is \citet{gnettner:kirch:2025}, who study confidence
sequences for weakly dependent time series. However, their approach
relies on long-run variance estimation, is restricted to the mean
parameter, and uses boundary assumptions that exclude the optimal
law-of-the-iterated-logarithm (LIL) width rate. Another recent extension was developed by \citet{ge:liu:2026} who propose confidence sequences for
means of dependent, potentially nonstationary time series. Anytime-valid inference
for nonlinear parameters under general weak dependence remains
largely undeveloped.\\
Another related research direction develops sequential estimation
procedures with prescribed precision-based stopping criteria, as in
\citet{flegal:gong:2015} and the Sequest procedure of
\citet{alexopoulos:goldsman:mokashi:tien:wilson:2019}.
The stopping criterion is baked into the procedure; in particular, the
asymptotic coverage guarantees of \citet{flegal:gong:2015} concern their confidence intervals at the
specified stopping time. Confidence sequences allow the stopping
strategy to be revised as observations accumulate, while preserving
the simultaneous coverage guarantee. For example, after obtaining
evidence that an effect is positive, a researcher may decide to
continue collecting data to establish that the same effect exceeds a
practically relevant threshold and, encouraged by subsequent evidence,
pursue a still larger threshold. If this last goal appears too
demanding, the researcher can stop and report the conclusions already
supported by the confidence sequence. Such adjustments to research goals as more data are observed are common
in scientific practice.\\
A key ingredient for our proposed methods is self-normalization
for dependent time series. Early contributions in this direction
include \citet{kiefer:vogelsang:bunzel:2000}, who developed tests for
regression hypotheses under dependence, and \citet{lobato:2001}, who
proposed a self-normalized test for uncorrelatedness under dependence.
Self-normalized methods are attractive because they avoid the bandwidth
choices involved in conventional long-run variance estimation and
dependent bootstrap procedures, often provide improved finite-sample
size and coverage accuracy \citep{lobato:2001,shao:2015}, and can
improve robustness to perturbations of the autocorrelation structure \citep{mueller:2007}.
Later, \citet{shao:2010} extended self-normalized inference to general
approximately linear estimators. A helpful overview of
self-normalization is given in \citet{shao:2015}.\\
Motivated by the benefits of self-normalization for retrospective inference, we extend it to anytime-valid inference in this paper, covering a broad class
of time-series parameters. The main theoretical challenge is controlling the fluctuations of the parameter estimator, the normalizer and the linearization error over an infinite time horizon.
Our main contributions are as follows:
\begin{itemize}
\item \textbf{Asymptotically sharp confidence sequences under weak dependence.}
We establish conditions on dependence and estimator structure under
which self-normalized confidence sequences have simultaneous coverage
over all times $n\ge m$ converging to $1-\alpha$ as $m\to\infty$.
The construction avoids estimation of the long-run variance and the
associated parameter tuning.

\item \textbf{General parameters.}
Our theory covers large classes of parameters that are functionals of
the marginal distribution, including quantiles, expectiles and
expected shortfall. Beyond marginal features, our theory also covers
lagged parameters, as we illustrate with the lag-one linear prediction coefficient.

\item \textbf{Two sequential self-normalizers.}
Motivated by existing retrospective normalizers, we propose two
natural extensions to the sequential setting, based on different
interpretations of what a self-normalizer does. First, we discuss the forward normalizer, which
uses only the history of sequential parameter estimates, allowing
inference without computing additional block estimators. This may be most useful where additional estimates are either expensive or restricted (we discuss an example involving differential privacy).
Second, we discuss the
three-block normalizer, which draws on additional block estimates. The three-block variant typically delivers narrower confidence sequences and it is compatible with the optimal LIL
width rate of intervals.
\end{itemize}
\noindent
We investigate the finite-sample coverage and interval widths of
both constructions through simulations. Mathematical proofs and an application of our method to Markov chain Monte Carlo output analysis are given in the supplementary appendix.
The appendix is not included in this arXiv version.

\section{Statistical Methodology}
\label{sec:method}
In this section, we develop self-normalized asymptotic confidence sequences for stationary time series. We begin with the notation and the general self-normalized construction principle, before introducing two sequential self-normalizers in separate subsections.\\
Let $(X_t)_{t\in\mathbb N}$ be a strictly stationary time series with marginal cumulative distribution function (cdf) $F$, and let $\theta\in\mathbb R$ be the parameter of interest. For integers $1\le a\le b$, write $a{:}b=\{a,\ldots,b\}$ and let $\widehat\theta_{a:b}$ denote an estimator based on the observations $X_a,\ldots,X_b$. A leading class of examples consists of marginal parameters $\theta=T(F)$, where $T$ is a potentially nonlinear functional. In this case, we consider the empirical cdf and the corresponding plug-in estimator
\begin{align}\label{e:emp:cdf}
\widehat F_{a:b}(x)
=\frac{1}{b-a+1}\sum_{i=a}^b\mathbf1\{X_i\le x\},
\qquad
\widehat\theta_{a:b}=T(\widehat F_{a:b}).
\end{align}
Examples include quantiles, expected shortfall, and expectiles (with the mean as a special case). The general notation $\widehat\theta_{a:b}$ also accommodates parameters involving serial relationships, such as autoregressive coefficients.
Our objective is to construct asymptotic confidence sequences in the sense of \eqref{eq:intro-asymptotic-cs}. A starting point for both constructions is an asymptotic linearization
\begin{equation}\label{eq:method-linearization}
\widehat\theta_{a:b}-\theta
=\frac{1}{b-a+1}\sum_{t=a}^b Y_t+R_{a:b},
\end{equation}
where $(Y_t)_{t\in\mathbb Z}$ is a centered stationary influence process and $R_{a:b}$ is a remainder. For marginal parameters, we typically have $Y_t=\psi_T(X_t)$, where $\psi_T$ is the influence function of $T$ at $F$. Under suitable conditions on dependence and the remainders, the leading partial sums admit a Brownian approximation, which we write informally as
\begin{align}\label{e:heurW}
n(\widehat\theta_{1:n}-\theta)
\approx \sum_{t=1}^n Y_t
\approx \sigma W(n),
\qquad n\ge m,
\end{align}
where $(W(t))_{t\ge0}$ is a standard Brownian motion and
\begin{align}\label{e:LRV}
\sigma^2
=\operatorname{Var}(Y_1)
+2\sum_{h=1}^{\infty}
\operatorname{Cov}(Y_1,Y_{1+h})
\end{align}
is the long-run variance of the influence process. The required control of these approximations over blocks and over the infinite time horizon will be mathematically studied in Section~\ref{sec:theory}.\\
One approach to sequential inference is to estimate $\sigma$ and use the Brownian approximation \eqref{e:heurW} to calibrate the intervals, as done for the mean parameter by \citet{gnettner:kirch:2025}. For nonlinear parameters, this can require estimating the influence function, which depends on the unknown marginal distribution, as well as accounting for serial dependence, which is difficult.

As an alternative we suggest using self-normalization, which we now roughly explain (more details are given in the following sections).  The idea is to use differences between estimates from selected blocks $a{:}b\subseteq1{:}n$ to construct a nonnegative normalizer $\widehat V_{n,m}$. The essential property is that the parameter estimator $\widehat\theta_{1:n}$ and the normalizer $\widehat V_{n,m}$ have fluctuations proportional to $\sigma$. Importantly, $\widehat V_{n,m}$ does not estimate $\sigma$ and, for any bounded ratio $t=n/m$, does not converge to a fixed number. Rather, it has a non-degenerate random limit $\sigma\cdot V(t)$, where the distribution of $V(t)$ is free of unknown model parameters (it is a functional of the standard Brownian motion).
Using Brownian scaling with $t=n/m$, and writing $W$ again for the rescaled Brownian motion, we obtain informally
\begin{equation}\label{eq:method-joint-approximation}
\frac{\sqrt n\,(\widehat\theta_{1:n}-\theta)}{\widehat V_{n,m}}
\approx
\frac{\cancel{\sigma}\,W(t)}
{\cancel{\sigma}\,\sqrt t\,V(t)},
\qquad t\ge1.
\end{equation}
The above approximation serves as the basis for confidence intervals.
Let $\rho:[1,\infty)\to(0,\infty)$ be a deterministic weight function, which determines how the usual $n^{-1/2}$ width is adjusted for simultaneous coverage over time. Now, choose $\rho$ so that
\begin{equation}\label{eq:method-brownian-supremum}
S_{\rho,V}
:=\sup_{t\ge1}
\frac{|W(t)|}{\rho(t)\, \sqrt t\,V(t)}
\end{equation}
is almost surely finite, and let $q_{1-\alpha}^{V}$ be its $(1-\alpha)$-quantile. Crucially, $S_{\rho,V}$ (and thus $q_{1-\alpha}^{V}$) is free of unknown model parameters and can be readily simulated. We then define the confidence sequence
\begin{equation}\label{eq:method-self-normalized-cs}
C_{n,m}
:=
\left[
\widehat\theta_{1:n}
-q_{1-\alpha}^{V}
\frac{\rho(n/m)\widehat V_{n,m}}{\sqrt n},
\quad
\widehat\theta_{1:n}
+q_{1-\alpha}^{V}
\frac{\rho(n/m)\widehat V_{n,m}}{\sqrt n}
\right],
\qquad n\ge m.
\end{equation}
The theory in Section~\ref{sec:theory} gives conditions on dependence, estimator linearization, and the weight function under which the simultaneous coverage probability converges to $1-\alpha$ as $m\to\infty$.

\noindent The two constructions of $\widehat V_{n,m}$ that we present in this paper use self-normalization in different ways. In Section~\ref{sec:forward}, we introduce the forward self-normalizer, starting from retrospective self-normalization for the mean. This is the most immediate sequential extension and requires only the history of sequential parameter estimates. Its simplicity comes with a modest rate limitation: ensuring that $S_{\rho,V}<\infty$, as required for our theory, here excludes the optimal law-of-the-iterated-logarithm (LIL) rate function (defined in \eqref{eq:three-block-lil-weight}), even though close to optimal weights are allowed (see \eqref{e:WF}). This means that the width rate for this normalizer is slightly suboptimal. 
To achieve the optimal width rate, Section~\ref{sec:3block} introduces the three-block self-normalizer. We motivate this construction by revisiting retrospective self-normalization through its fixed-$\kappa$ interpretation, which suggests a different, more surprising way to use the available subsamples.

\subsection{The forward self-normalizer}
\label{sec:forward}

We begin with a short review of self-normalization in retrospective inference, which motivates our first sequential construction.

\noindent\textbf{Motivation part I.}
Consider observations $X_1,\ldots,X_m$ from a weakly dependent, stationary time series. The aim is to construct a confidence interval for the mean parameter $\theta:=\mathbb E X_1$. The corresponding block estimators are
\[
\widehat\theta_{a:b}
:=\bar X_{a:b}
:=\frac{1}{b-a+1}\sum_{t=a}^b X_t,
\]
and the forward empirical means $\bar X_{1:1},\ldots,\bar X_{1:m}$ are already sufficient for self-normalization. Specifically, the classical self-normalizer is
\begin{equation}\label{eq:retrospective-mean-sn}
\widehat V_m
:=
\left\{
\frac{1}{m^2}\sum_{k=1}^m k^2
\bigl(\bar X_{1:k}-\bar X_{1:m}\bigr)^2
\right\}^{1/2};
\end{equation}
see \citet{lobato:2001} and \citet{shao:2010}.
Under the usual functional central limit theorem for the partial sums, the scaled deviations
$k(\bar X_{1:k}-\bar X_{1:m})/\sqrt m$
converge to a Brownian bridge on the time scale $k/m$. More precisely, as $m\to\infty$, we have the joint convergence
\begin{equation}\label{e:mot:1}
\left(
\sqrt m\,(\bar X_{1:m}-\theta),\,
\widehat V_m
\right)
\ \Rightarrow\
\left(
\sigma W(1),\,
\sigma
\left\{\int_0^1
\bigl(W(r)-rW(1)\bigr)^2\,dr
\right\}^{1/2}
\right),
\end{equation}
where $W$ is a standard Brownian motion and $\sigma^2>0$ is the long-run variance of $(X_t)$. Thus, the unknown scale cancels in the self-normalized statistic:
\begin{equation}\label{eq:retrospective-mean-pivot}
\frac{\sqrt m\,(\bar X_{1:m}-\theta)}{\widehat V_m}
\ \Rightarrow\
\frac{W(1)}
{\left\{\int_0^1
\bigl(W(r)-rW(1)\bigr)^2\,dr
\right\}^{1/2}}
=:Z_{\mathrm{SN}}.
\end{equation}
The limiting distribution is non-normal but pivotal. Its quantiles therefore allow hypothesis testing and confidence interval construction.

\smallskip
\noindent\textbf{The forward normalizer.}
A straightforward extension to the sequential setting replaces the empirical means in \eqref{eq:retrospective-mean-sn} by general parameter estimates and evaluates the normalizer at each time $n$. We thus define
\begin{equation}\label{eq:forward-sn}
\widehat V_n^{\mathrm{fw}}
:=
\left\{
\frac{1}{n^2}\sum_{k=1}^n k^2
\bigl(\widehat\theta_{1:k}-\widehat\theta_{1:n}\bigr)^2
\right\}^{1/2}.
\end{equation}
Its Brownian counterpart on the rescaled time axis is
\begin{equation}\label{eq:forward-brownian-normalizer}
V^{\mathrm{fw}}(x)
:=
\frac1x
\left\{
\int_0^x
\left(W(u)-\frac uxW(x)\right)^2\,du
\right\}^{1/2},
\qquad x\ge1.
\end{equation}
In Section~\ref{sec:theory}, we establish conditions on dependence, estimator linearization, and the growth of $\rho$ under which
\begin{equation}\label{e:lim:1}
\sup_{n\ge m}
\frac{\sqrt n\,|\widehat\theta_{1:n}-\theta|}
{\rho(n/m)\widehat V_n^{\mathrm{fw}}}
\ \Rightarrow\
\sup_{x\ge1}
\frac{|W(x)|}
{\sqrt x\,\rho(x)V^{\mathrm{fw}}(x)}
\qquad\text{as }m\to\infty.
\end{equation}
Let $q_{1-\alpha}^{\mathrm{fw}}$ denote the $(1-\alpha)$-quantile of the limiting distribution, with its dependence on $\rho$ suppressed in the notation. Substituting this critical value and $\widehat V_n^{\mathrm{fw}}$ into \eqref{eq:method-self-normalized-cs} yields the forward confidence sequence. Under the conditions given in Section~\ref{sec:theory}, its simultaneous coverage probability converges to $1-\alpha$.

\smallskip
\noindent\textbf{Computational properties.}
The main computational advantage of \eqref{eq:forward-sn} is that it uses only the sequential parameter estimates $\{\widehat\theta_{1:n}\}_{n\in\mathbb N}$. No additional block estimators need to be computed. Moreover, the normalizer admits a simple recursive update, namely with
\[
A_n:=A_{n-1}+n^2\widehat\theta_{1:n},
\qquad
B_n:=B_{n-1}+n^2\widehat\theta_{1:n}^2,
\qquad A_0=B_0=0,
\]
we have
\begin{equation}\label{eq:forward-update}
\bigl(\widehat V_n^{\mathrm{fw}}\bigr)^2
=
\frac{
B_n-2\widehat\theta_{1:n}A_n
+\widehat\theta_{1:n}^2\,n(n+1)(2n+1)/6
}{n^2}.
\end{equation}
Thus, once the next estimate is available, updating the normalizer requires only constant additional computation and storage; it is not even necessary to save the complete estimator history.
The advantages of the forward normalizer are particularly relevant when additional block estimates are expensive to compute, when access to the underlying observations is restricted, or when further releases of data-dependent quantities are undesirable. The following example illustrates the last situation.

\begin{example}[Inference from differentially private estimates]
\label{ex:forward-dp}
Aggregate statistics can unwittingly reveal information about the observations
that gave rise to them -- even when the individual data points were never
published. This concern has given rise to the concept of differential
privacy (DP), which is widely used to control information leakage
through the publication of statistical results. Typically, DP is
achieved by adding suitably calibrated random noise to statistics
before publication; see \citet{slavkovic:seeman:2023} for a statistical
review. Forward self-normalization is particularly useful when
inference must be based on such privatized estimates.
Suppose an analyst receives a stream of sequentially updated, but
privatized, parameter estimates
$\{\widehat\theta_{1:n}^{\mathrm{DP}}\}_{n\in\mathbb N}$,
and that the entire published sequence satisfies
$\varepsilon$-differential privacy (with respect to replacement of a
single observation). Roughly, this means that changing one
observation $X_i$ has only a limited effect on the distribution of the
published results, with $\varepsilon$ controlling the strength of
this protection.\\
A private, sequential publication scheme might look as follows:  Construct recursively (starting from a fixed starting value)
\[
\widehat\theta_{1:n}^{\mathrm{DP}}
=M_n\bigl(\widehat\theta_{1:n-1}^{\mathrm{DP}},X_n\bigr).
\]
Here, $M_n$ is a randomized mechanism that updates $\widehat\theta_{1:n-1}^{\mathrm{DP}}$ using the current value $X_n$, while maintaining $\varepsilon$-DP for $X_n$.  Since $X_n$ is
only accessed once directly, later releases of parameter estimates $\widehat\theta_{1:{n+k}}^{\mathrm{DP}}$ do not reveal additional information about $X_n$ (so-called "post-processing stability"; Proposition~2.1 in \citet{dwork:roth:2014}).
The forward normalizer can be computed entirely from the published
parameter estimates
$\{\widehat\theta_{1:n}^{\mathrm{DP}}\}_{n\in\mathbb N}$.
It does not require any additional direct access to the  observations, and again due to post-processing stability, computing and publishing normalizers maintains the $\varepsilon$-DP guarantee.
This means that the confidence intervals from the forward normalizer provide a "free lunch" in the sense that they provide statistical inference without any additional privacy costs.\\
In contrast, computing and publishing additional parameter estimates
(from more subsets of the data) risks leaking  extra information. This extra leakage is a major obstacle to alternatives to self-normalization, such as nonparametric
bootstrap under DP. For a discussion of these difficulties  in the retrospective setting
see \citet{dette:graw:2025}.
\end{example}

\smallskip
\noindent\textbf{Mathematical properties.}
The normalizer $\widehat V_n^{\mathrm{fw}}$ does not depend on the burn-in period $m$. Nevertheless, the resulting confidence sequence does depend on $m$ through its starting time and the rescaling in the weight function $\rho(n/m)$. A sufficiently large burn-in is generally needed for the Brownian calibration to be accurate, as formalized by \eqref{e:lim:1}. Special models permit exact finite-sample calibration. For the mean of i.i.d. Gaussian observations, for example, the discrete self-normalized process is pivotal for each fixed $m\ge2$. Its finite-sample critical value can therefore be simulated directly. Our asymptotic theory goes far beyond this case, allowing general data distributions, nonlinear parameters and  serial dependence.\\
The weight function $\rho$ supplies the deterministic growth factor in the interval half-width
\[
q_{1-\alpha}^{\mathrm{fw}}(\rho)\,
\rho(n/m)\,
\frac{\widehat V_n^{\mathrm{fw}}}{\sqrt n}.
\]
Among admissible choices, slower growth of $\rho$ improves the width rate as $n/m\to\infty$. If the true long-run standard deviation $\sigma$ were available, the law of the iterated logarithm would identify $\sqrt{\log\log x}$ as the smallest admissible growth order for $\rho(x)$ at large $x$. The forward self-normalizer introduces additional fluctuations. At any fixed $x$, $V^{\mathrm{fw}}(x)$ has the same nondegenerate distribution as the retrospective Brownian normalizer. Over an unbounded horizon, however, it can become arbitrarily small. Time-uniform calibration must account for these small-denominator excursions as well as the fluctuations of the numerator.\\
Theorem~\ref{thm:main} shows that the weight function
\begin{align}\label{e:WF}
\rho_{\mathrm{fw}}(x):=\log(\mathrm e+\log x),\qquad x\ge1,
\end{align}
is sufficient. Its log-log growth is slightly larger than the square-root log-log growth of a LIL boundary. This raises a natural question: can a sequential self-normalizer cancel the unknown long-run variance while attaining the optimal LIL width rate? The next construction achieves precisely this.

\subsection{The three-block self-normalizer}
\label{sec:3block}

The three-block self-normalizer uses blocks whose length is tied to the
burn-in period. As monitoring continues, the normalizer can therefore
average over increasingly many blocks of that length. This leads to
greater stability over long horizons and allows confidence sequences
with the LIL width rate. To motivate the construction, we revisit
retrospective self-normalization from a different perspective.

\noindent\textbf{Motivation part II.}
We return to inference for the mean parameter $\theta=\mathbb E X_1$
and give a second interpretation of the retrospective normalizer in
\eqref{eq:retrospective-mean-sn}. For observations $X_1,\ldots,X_m$,
define the sample autocovariances
\[
\widehat\gamma_m(h)
:=
\frac1m\sum_{t=1}^{m-h}
(X_t-\bar X_{1:m})(X_{t+h}-\bar X_{1:m}),
\qquad h=0,\ldots,m-1.
\]
The Bartlett-kernel estimator of the long-run variance with
bandwidth $\ell\in\{1,\ldots,m\}$ is
\begin{equation}\label{eq:bartlett-lrv}
\widehat\sigma_m^2(\ell)
:=
\widehat\gamma_m(0)
+2\sum_{h=1}^{\ell-1}
\left(1-\frac h\ell\right)\widehat\gamma_m(h).
\end{equation}
Under suitable conditions, this estimator is consistent when
$\ell\to\infty$ and $\ell/m\to0$. A conventional bandwidth choice
is of order $m^{1/3}$, which balances the leading squared bias and
variance terms \citep{andrews:1991}. An interesting case arises when the
bandwidth is kept proportional to the sample size, so that
$\ell/m\to \kappa\in(0,1]$. The resulting variance estimator becomes inconsistent, because it
has a nondegenerate limiting distribution. Nevertheless, its limiting distribution remains
proportional to the true long-run variance, so it becomes possible to cancel it out. 
Indeed, the development of self-normalization for time series is
rooted in the early work on inference with inconsistent long-run
variance estimators \citep{kiefer:vogelsang:bunzel:2000,lobato:2001},
which also laid the foundations for the fixed-$\kappa$ framework of
\citet{kiefer:vogelsang:2005}. For the mean parameter, the connection
is particularly transparent at bandwidth $\ell=m$. As shown by
\citet{kiefer:vogelsang:2002}, the Bartlett estimator then satisfies
the exact identity
\begin{equation}\label{eq:bartlett-retrospective-sn}
\frac12\widehat\sigma_m^2(m)
=
\frac1{m^2}\sum_{k=1}^m
k^2\bigl(\bar X_{1:k}-\bar X_{1:m}\bigr)^2
=
\widehat V_m^2.
\end{equation}
Thus, up to the factor $1/2$, the squared retrospective
self-normalizer is nothing but a Bartlett estimator with bandwidth $\ell=m$ ($\kappa=1$). Its joint limit with the sample mean is given in
\eqref{e:mot:1}. For this exposition we fix $\kappa=1$, but we discuss other values of $\kappa$ at the end of this section.\\
The representation \eqref{eq:bartlett-retrospective-sn} suggests two natural sequential extensions.
At time $n\ge m$, the forward normalizer uses bandwidth $\ell=n$
and hence keeps the bandwidth-to-sample-size ratio equal to one.
The second construction sticks with the initial bandwidth $\ell=m $,
while using all $n$ observations to compute the sample
autocovariances.  
The bandwidth $\ell$ then remains fixed as monitoring
proceeds, while its ratio to the current sample size decreases with 
$m/n$. For the mean, this gives the normalizer
\begin{equation}\label{eq:three-block-bartlett}
\widehat V_{n,m}^{\mathrm{3b}}
:=
\left\{\frac12\widehat\sigma_n^2(m)\right\}^{1/2},
\qquad n\ge m,
\end{equation}
where $\widehat\sigma_n^2(m)$ is defined by 
\eqref{eq:bartlett-lrv}.
Expanding the Bartlett quadratic form and grouping its terms gives
the following  three-block representation
\begin{equation}\label{eq:three-block-mean-sn}
\begin{split}
\bigl(\widehat V_{n,m}^{\mathrm{3b}}\bigr)^2
=\frac1{2nm}\bigg\{&
\sum_{k=1}^{m-1} k^2
\bigl(\bar X_{1:k}-\bar X_{1:n}\bigr)^2+\sum_{j=0}^{n-m}m^2
\bigl(\bar X_{j+1:j+m}-\bar X_{1:n}\bigr)^2\\
&+\sum_{k=1}^{m-1}k^2
\bigl(\bar X_{n-k+1:n}-\bar X_{1:n}\bigr)^2
\bigg\}.
\end{split}
\end{equation}
The three blocks explain the name \textit{three-block self-normalizer}.
At $n=m$, the middle one vanishes and the two boundary sums
coincide. Hence, for the mean, and at $n=m$
\[
\widehat V_{m,m}^{\mathrm{3b}}
=\widehat V_m
=\widehat V_m^{\mathrm{fw}}.
\]
Both sequential constructions thus recover the same retrospective
normalizer at the end of the burn-in period, while their subsequent
behavior differs critically.

\medskip
\noindent\textbf{The three-block normalizer.}
For general parameters, we replace the block means in
\eqref{eq:three-block-mean-sn} by the corresponding parameter
estimates and define, for $n\ge m\ge2$,
\begin{equation}\label{eq:three-block-sn}
\begin{split}
\bigl(\widehat V_{n,m}^{\mathrm{3b}}\bigr)^2
:=\frac1{2nm}\bigg\{&
\sum_{k=1}^{m-1} k^2
\bigl(\widehat\theta_{1:k}-\widehat\theta_{1:n}\bigr)^2+\sum_{j=0}^{n-m}m^2
\bigl(\widehat\theta_{j+1:j+m}
-\widehat\theta_{1:n}\bigr)^2\\
&+\sum_{k=1}^{m-1}k^2
\bigl(\widehat\theta_{n-k+1:n}
-\widehat\theta_{1:n}\bigr)^2
\bigg\}.
\end{split}
\end{equation}
The moving-block contribution (middle sum) is central to the long-horizon
behavior of this construction. Each new observation adds another
block of length $m$, so the normalizer averages over increasingly
many blocks as $n/m$ grows. Of course, these blocks overlap and are dependent,
but the averaging still reduces their variance as $n/m \to \infty$. The two boundary
contributions become negligible at late times. So, the entire construction stabilizes for very late values $n/m \to \infty$.

\medskip
\noindent\textbf{Sequential calibration.}
Let $W$ be a standard Brownian motion and write
\[
G_x(u):=W(u)-\frac uxW(x),
\qquad 0\le u\le x,\quad x\ge1.
\]
The Brownian counterpart of \eqref{eq:three-block-sn} is
\begin{equation}\label{eq:three-block-brownian-normalizer}
\begin{split}
\bigl(V^{\mathrm{3b}}(x)\bigr)^2
:=\frac1{2x}\bigg\{&
\int_0^1 G_x(u)^2\,du+\int_0^{x-1}
\bigl(G_x(u+1)-G_x(u)\bigr)^2\,du
+\int_{x-1}^x G_x(u)^2\,du
\bigg\}.
\end{split}
\end{equation}
Under the dependence, linearization, and weight conditions in
Section~\ref{sec:theory}, we have
\begin{equation}\label{eq:three-block-sequential-limit}
\sup_{n\ge m}
\frac{\sqrt n\,|\widehat\theta_{1:n}-\theta|}
{\rho(n/m)\widehat V_{n,m}^{\mathrm{3b}}}
\ \Rightarrow\
\sup_{x\ge1}
\frac{|W(x)|}
{\sqrt x\,\rho(x)V^{\mathrm{3b}}(x)}
\qquad\text{as }m\to\infty.
\end{equation}
Let $q_{1-\alpha}^{\mathrm{3b}}$ denote the $(1-\alpha)$-quantile
of the limiting distribution. Then, substituting $q_{1-\alpha}^{\mathrm{3b}}$
and $\widehat V_{n,m}^{\mathrm{3b}}$ into
\eqref{eq:method-self-normalized-cs} gives an asymptotic $(1-\alpha)$ confidence sequence.

\medskip
\noindent\textbf{Long-horizon behavior and interval width.}
The main statistical advantage of the three-block construction is
captured by the stabilization of its Brownian counterpart:
\begin{equation}\label{eq:three-block-stabilization}
\bigl(V^{\mathrm{3b}}(x)\bigr)^2
\longrightarrow\frac12
\qquad\text{almost surely as }x\to\infty.
\end{equation}
Thus, $V^{\mathrm{3b}}(x)$ stays bounded away from zero at sufficiently large rescaled times - unlike the forward normalizer. This stability permits more slowly growing weights $\rho$ and thereby improves the achievable interval-width at late times $n$.
In particular, our theory permits the optimal LIL-type weight
\begin{equation}\label{eq:three-block-lil-weight}
\rho_{\mathrm{LIL}}(x)
:=2\sqrt{\log(\mathrm e+\log x)},
\qquad x\ge1,
\end{equation}
with corresponding rate of the interval width 
\[
\frac{\rho_{\mathrm{LIL}}(n/m)}{\sqrt n}
\asymp
\sqrt{\frac{\log\log(n/m)}{n}},
\qquad n/m\to\infty.
\]
This improves on the boundary growth required by the forward
normalizer. In finite samples, we find that the three-block normalizer also has better coverage properties for small to moderate burn-in periods $m$.

\begin{figure}[htbp]
\centering
\includegraphics[width=0.65\textwidth]{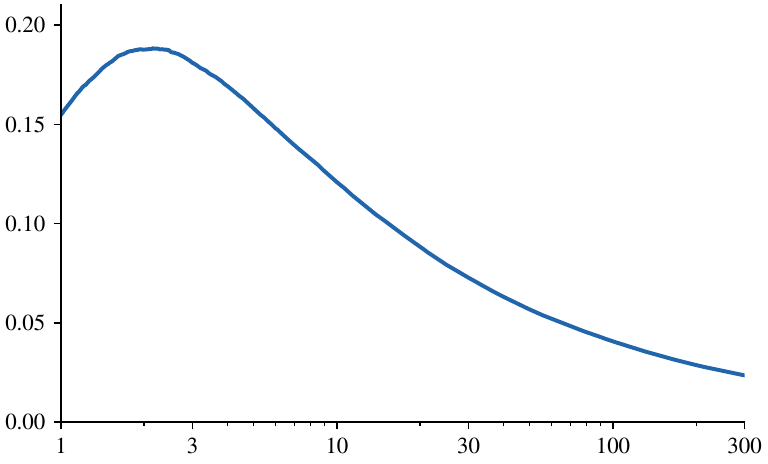}
\caption{Standard deviation of the three-block Brownian normalizer
$V^{\mathrm{3b}}(x)$ for $x\in[1,300]$, estimated from
$50{,}000$ Brownian paths with $200$ increments per unit of time.
The time axis is logarithmic. The normalizer remains random
at every fixed $x$, while its variability eventually decreases as $x$ increases.}
\label{fig:three-block-sd}
\end{figure}

\medskip
\noindent\textbf{Are we just estimating the long-run variance through the back door?} At first glance, the distinction between three-block self-normalizer and long-run variance estimation might seem to be merely a matter of bandwidth selection: for the mean parameter, the self-normalizer has the Bartlett representation \eqref{eq:three-block-bartlett}, with bandwidth $\ell=m$.\footnote{Of course, for estimation, this would be  a rather odd bandwidth choice. For estimation with $m=n=200$, a more standard bandwidth would be e.g. $\ell=\lfloor m^{1/3} \rfloor=5$ rather than the extreme $\ell=m=200$.} This impression can be reinforced when we notice that, for large $m$, and $n=\lfloor xm\rfloor$ 
\begin{align} \label{e:conheur}
 \big(\widehat V_{\lfloor mx\rfloor,m}^{\mathrm{3b}}\big)^2  \approx \sigma^2 \bigl(V^{\mathrm{3b}}(x)\bigr)^2
\longrightarrow\frac{\sigma^2}2
\qquad\text{ as }x\to\infty.
\end{align}
Here $m\to\infty$ is taken first, and the subsequent convergence as $x\to\infty$ follows from \eqref{eq:three-block-stabilization}. Does this not mean that the normalizer eventually estimates the long-run variance? Only in the limited sense of \eqref{e:conheur}, and only for very late times (here, large $x$). But our confidence
sequence does not wait for this stabilization to happen. Inference already
begins for small ratios $x \ge 1$, i.e. when the normalizer still has substantial random fluctuations,
and our calibration explicitly accounts for these fluctuations.  
Figure~\ref{fig:three-block-sd} illustrates this. It shows that while variability of $V^{\mathrm{3b}}(x)$ eventually decays, it remains non-negligible,  even for quite large ratios of $x = n/m$.
Accurate long-run variance estimation would require this variability
to be uniformly small - self-normalized inference remains valid throughout this
transition. \\
Even though stabilization occurs only late, it still has important implications for our mathematical theory. In particular, it makes the three-block normalization compatible with the rate-optimal LIL weights, that were not compatible with the earlier discussed forward normalizer. This translates into narrower confidence intervals for late time points.

\medskip
\noindent\textbf{Computational properties.}
The three-block normalizer can be updated sequentially.
For this purpose, suppose we store the
running sums
\[
A_{n,m}
:=\sum_{j=0}^{n-m}\widehat\theta_{j+1:j+m},
\qquad
B_{n,m}
:=\sum_{j=0}^{n-m}\widehat\theta_{j+1:j+m}^{\,2}.
\]
The middle block in \eqref{eq:three-block-sn} is then
\[
m^2\left\{
B_{n,m}-2\widehat\theta_{1:n} A_{n,m}
+(n-m+1)\widehat\theta_{1:n}^2
\right\}.
\]
With the next observation, only the new block estimate needs to be
added:
\[
A_{n+1,m}
=A_{n,m}+\widehat\theta_{n-m+2:n+1},\qquad 
B_{n+1,m}
=B_{n,m}+\widehat\theta_{n-m+2:n+1}^{\,2}.
\]
Thus, earlier middle-block estimates need not be recomputed or
stored individually. The first block contribution can
similarly be evaluated from weighted sums computed once during
the burn-in period. The last block contribution uses the current
$m-1$ partial-block estimates. For fixed $m$, the number of
additional block estimates needed at each update therefore does
not grow with the monitoring horizon.

\begin{figure}[htbp]
\centering
\includegraphics[width=0.65\textwidth]{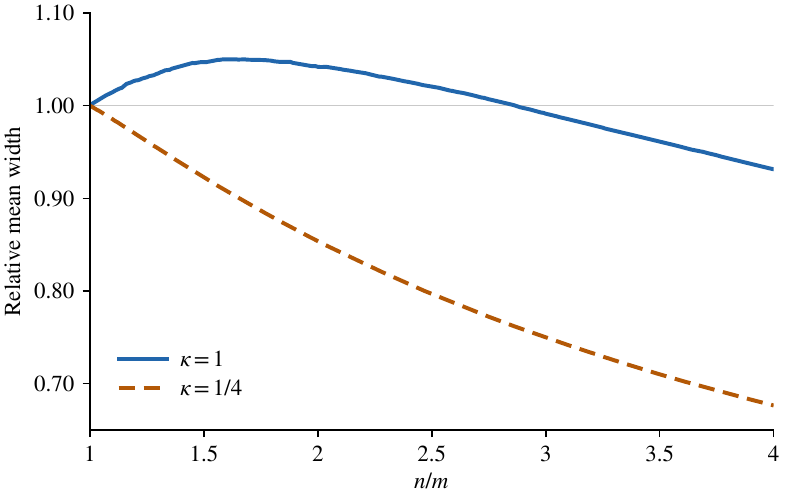}
\caption{Relative mean interval width for the three-block method with
the LIL weight, for different $\kappa$ values. 
For each $\kappa$, we divide the mean width at time $n$ by the mean width at $n=m$ (so both curves start at $1$). Widths are simulated using the Brownian limit and on the $x$-axis we report rescaled time $t=n/m$. Both curves are based on $50{,}000$
Brownian paths with $200$ increments per unit of time. \label{fig:three-block-relative-mean-width}}

\end{figure}

\textbf{Fixed-$\kappa$ asymptotics.}
In the above motivation of the three-block method, we fix the
bandwidth at $\ell=m$. Following the general logic of fixed-$\kappa$
asymptotics, we may instead choose
$\ell=\lfloor\kappa m\rfloor$ for any fixed $\kappa\in(0,1]$.
More precisely, define
\[
\widehat V_{n,m,\kappa}^{\mathrm{3b}}
:=\widehat V_{n,\ell}^{\mathrm{3b}},
\qquad n\ge m,
\]
where the right-hand side is given by
\eqref{eq:three-block-sn} with every occurrence of $m$ replaced
by $\ell$. Monitoring still begins at $n=m$.
Under Assumption~\ref{ass:general} and the three-block weight
condition in \eqref{eq:admissible-weights}, the same argument gives
\[
\sup_{n\ge m}
\frac{\sqrt n\,|\widehat\theta_{1:n}-\theta|}
{\rho(n/m)\widehat V_{n,m,\kappa}^{\mathrm{3b}}}
\ \Rightarrow\
\sup_{x\ge1/\kappa}
\frac{|W(x)|}
{\sqrt x\,\rho(\kappa x)V^{\mathrm{3b}}(x)}
\qquad\text{as }m\to\infty.
\]
Here $V^{\mathrm{3b}}$ is the Brownian normalizer from
\eqref{eq:three-block-brownian-normalizer}, driven by the same
Brownian motion $W$ as the numerator.
The critical values must, of course, be calibrated accordingly. Why would one consider values other than the (generic sounding) $\kappa=1$?
The answer is that at $\kappa=1$ the confidence sequence can have some counterintuitive evolution: confidence intervals can initially be narrow (at
$n\approx m$), then widen as the normalizer stabilizes, before
eventually shrinking again for larger $n$. This behavior seems counterintuitive
because, at least temporarily, more data is associated with wider intervals. For $\kappa=1$, the expected normalizer initially grows fast enough for the mean interval width to increase (see the blue curve in Figure~\ref{fig:three-block-relative-mean-width}).
For $\kappa=1/4$, the mean width decreases throughout the displayed range (see the dashed curve in Figure~\ref{fig:three-block-relative-mean-width}).
One can view this choice of $\kappa$ as reallocating $\alpha$-spending in a more intuitive way.
This is similar
to, and indeed interlocked with, choosing a weight function, which also influences $\alpha$-spending.
It should be distinguished from tuning for consistent variance
estimation: at fixed $n/m$, any fixed $\kappa>0$ still yields an
inconsistent variance estimator. Based on the comparisons in
the supplementary appendix, we recommend $\kappa=1/4$ as a practical
default to govern $\alpha$-spending in conjunction with the
LIL weight function from this section.

\section{Theory}
\label{sec:theory}

In this section, we establish asymptotic time-uniform coverage for
the self-normalized confidence intervals in
\eqref{eq:method-self-normalized-cs}, using the forward and three-block
normalizers defined in \eqref{eq:forward-sn} and
\eqref{eq:three-block-sn}, respectively. Our inferential goal
is the asymptotic confidence-sequence guarantee in
\eqref{eq:intro-asymptotic-cs}.
We begin with high-level theoretical results for an abstract parameter $\theta$, before turning to concrete examples in the later part of this section.

\subsection{General theory}

For our theory, let $(X_i)_{i\in\mathbb Z}$ be the observation process. The reader should also recall the influence process $(Y_i)_{i\in\mathbb Z}$ from the estimator linearization \eqref{eq:method-linearization}.
To model temporal dependence along the time series, we invoke physical
dependence from \citet{wu:2005:dependence} (see the next definition).
For the statement, we use the standard $L^{a}$-norm $\|U\|_a=(\mathbb E|U|^a)^{1/a}$ for $a\ge1$.

\begin{definition}[Physical dependence]\label{def:causal}
Let $a\ge1$, and let $(Z_i)_{i\in\mathbb Z}$ be a real-valued
process satisfying $\mathbb E|Z_0|^a<\infty$ and admitting the
causal representation
\[
Z_i=G(\varepsilon_i,\varepsilon_{i-1},\ldots),
\qquad i\in\mathbb Z,
\]
where $(\varepsilon_i)_{i\in\mathbb Z}$ is a sequence of i.i.d.
innovations taking values in a metric space $S$, and
$G:S^{\mathbb N}\to\mathbb R$ is measurable.
Let $\varepsilon_0^\star$ have the same distribution as
$\varepsilon_0$ and be independent of the entire innovation
sequence. For $i\ge1$, define
\[
Z_i^\star
=G(\varepsilon_i,\varepsilon_{i-1},\ldots,
\varepsilon_1,\varepsilon_0^\star,\varepsilon_{-1},\ldots),
\]
with $Z_0^\star=G(\varepsilon_0^\star,\varepsilon_{-1},\ldots)$ and $Z_i^\star=Z_i$ for $i<0$.
The $L^a$-physical dependence coefficients associated with
this representation are
\[
\delta_a^Z(i)
:=
\left\|Z_i-Z_i^\star\right\|_a,
\qquad i\ge0.
\]
\end{definition}
$ $\\
The coefficient $\delta_a^Z(i)$ measures the influence of an innovation $i$ steps in the past on the current observation. 
Conditions on the decay of these coefficients provide a
quantitative description of weak dependence. We now summarize our main assumptions.

\begin{assumption}
\label{ass:general}
There exist $r>4$, $\eta\in(1/r,1/2)$, constants $C<\infty$
and $\chi\in(0,1)$, and a real-valued stationary causal process
$(Y_i)_{i\in\mathbb Z}$ such that the following conditions hold.

\begin{itemize}
\item[(i)] 
With physical dependence coefficients given in
Definition~\ref{def:causal}, the process $(Y_i)_{i\in\mathbb Z}$
satisfies
\[
\mathbb E Y_0=0,
\qquad
\|Y_0\|_r<\infty,
\qquad
\delta_r^Y(j)\le C\chi^j,
\quad j\ge0.
\]
Its long-run variance $\sigma^2$, defined in \eqref{e:LRV}, is strictly positive.

\item[(ii)] 
For every pair of integers $a\le b$, the block estimator
$\widehat\theta_{a:b}$ is measurable and finite almost surely,
and admits the representation \eqref{eq:method-linearization}.
The remainder array $(R_{a:b})_{a\le b}$ is jointly stationary
under integer time shifts, and the remainders satisfy
\begin{equation}\label{eq:block-remainder-bound}
\sup_{a\in\mathbb Z}
\|R_{a:a+\ell-1}\|_r
\le C\ell^{-1/2-\eta},
\qquad \ell\ge1.
\end{equation}
\end{itemize}
\end{assumption}

\noindent The influence process captures the first-order fluctuations
of the block estimators. For  functionals of the marginal distribution of $X_i$, it typically
takes the form $Y_i=\psi_T(X_i)$, where $\psi_T$ is an influence
function. More general forms occur, e.g., in our application to autoregressive processes; see Section~\ref{sec:examples}. 
For now, we develop a high-level theory, and thus impose conditions directly on the moments and weak dependence of the $Y_i$. Yet, in the examples below, we show how these conditions follow from moment and dependence properties of the underlying time series $(X_i)_{i\in\mathbb Z}$, which seems like the most natural level at which to impose them.

We can already state the main convergence and coverage result. The two constructions require different growth rates for the weight function.
The proofs for both constructions are given in the supplementary appendix.
For notational uniformity, write
$\widehat V_{n,m}^{\mathrm{fw}}:=\widehat V_n^{\mathrm{fw}}$.
If a denominator is zero, a ratio with positive numerator is interpreted as $+\infty$, and $0/0$ as zero. This convention matches the coverage event for an interval of zero width.

\begin{theorem}[Asymptotic time-uniform self-normalized inference]
\label{thm:main}
Suppose Assumption~\ref{ass:general} holds. Fix
$j\in\{\mathrm{fw},\mathrm{3b}\}$, and let
$\rho:[1,\infty)\to(0,\infty)$ be continuous and satisfy
\begin{equation}\label{eq:admissible-weights}
\begin{aligned}
\sup_{x\ge1}\frac{\rho_{\mathrm{fw}}(x)}{\rho(x)}<\infty
&\qquad\text{if }j=\mathrm{fw},\\
\sup_{x\ge1}\frac{\rho_{\mathrm{LIL}}(x)}{\rho(x)}<\infty
&\qquad\text{if }j=\mathrm{3b}.
\end{aligned}
\end{equation}
Here, $\rho_{\mathrm{fw}}$ is defined in \eqref{e:WF} and $\rho_{\mathrm{LIL}}$ in \eqref{eq:three-block-lil-weight}.
Then, as $m\to\infty$,
\begin{equation}\label{eq:main-convergence}
M_m^j(\theta)
:=\sup_{n\ge m}
\frac{\sqrt n\,|\widehat\theta_{1:n}-\theta|}
{\rho(n/m)\widehat V_{n,m}^j}
\ \Rightarrow\
S_{\rho,j}
:=\sup_{x\ge1}
\frac{|W(x)|}{\sqrt x\,\rho(x)V^j(x)},
\end{equation}
where $V^{\mathrm{fw}}$ and $V^{\mathrm{3b}}$ are given in
\eqref{eq:forward-brownian-normalizer} and
\eqref{eq:three-block-brownian-normalizer}, respectively.
The limits are finite almost surely, and their distributions depend only on $\rho$ and the choice of normalizer. 
Furthermore, for $\alpha\in(0,1)$, define
\[
q_{1-\alpha}^j
:=\inf\{q:\mathbb P(S_{\rho,j}\le q)\ge1-\alpha\}.
\]
If the distribution function of $S_{\rho,j}$ is continuous at
$q_{1-\alpha}^j$, the intervals \eqref{eq:method-self-normalized-cs}, with normalizer $\widehat V_{n,m}^j$ and critical value $q_{1-\alpha}^j$, satisfy
\begin{equation}\label{eq:main-coverage}
\lim_{m\to\infty}
\mathbb P\!\left(\theta\in C_{n,m}^j
\text{ for every }n\ge m\right)=1-\alpha.
\end{equation}
In particular, they form an asymptotic $(1-\alpha)$-confidence sequence in the sense of \eqref{eq:intro-asymptotic-cs}.
\end{theorem}

\begin{corollary}[Coverage at data-dependent times]
\label{cor:stopping}
Suppose the coverage conditions of Theorem~\ref{thm:main} hold, and let
$\tau_m$ be a random variable with values in $\mathbb{N}_{\ge m} \cup \{\infty\}$. Then
\[
\liminf_{m\to\infty}
\mathbb P\!\left(\theta\in C_{\tau_m,m}^j\right)\ge1-\alpha.
\]
\end{corollary}

Relevant stopping rules suggested in the literature include those where the present confidence interval is shorter than a pre-specified length (which may depend on $m$), or when a confidence interval does not include a hypothesized parameter value $\theta_0 $ for the first time. The latter case is equivalent to a sequential test over an infinite time horizon, and requires the strong guarantees of a confidence sequence. Under $\theta=\theta_0$, the probability of never stopping converges to $1-\alpha$ as $m\to\infty$.

\subsection{Examples}

\label{sec:examples}
We now give concrete examples of parameters and block estimators
and verify Assumption~\ref{ass:general} under sufficient conditions
on the underlying time series.

\noindent\textbf{Quantiles and expected shortfall.} 
Fix $p\in(0,1)$. Define the quantile and upper expected shortfall by
\begin{equation}\label{eq:quantile-es-parameters}
q_p:=F^{-1}(p),
\qquad
e_p:=q_p+\frac{\mathbb E(X_0-q_p)_+}{1-p},
\end{equation}
where $z_+=\max(z,0)$ and $F^{-1}(p)=\inf\{x:F(x)\ge p\}$. Their block estimators are
\begin{equation}\label{eq:quantile-es-estimators}
\widehat q_{p,a:b}:=\widehat F_{a:b}^{-1}(p),
\qquad
\widehat e_{p,a:b}
:=\widehat q_{p,a:b}
+\frac1{(b-a+1)(1-p)}\sum_{i=a}^b(X_i-\widehat q_{p,a:b})_+.
\end{equation}
When $F$ has a density $f$ that is continuous and positive at $q_p$, the corresponding influence processes are
\begin{equation}\label{eq:quantile-es-influences}
Y_i^q:=\frac{p-\mathbf1\{X_i\le q_p\}}{f(q_p)},
\qquad
Y_i^e:=q_p-e_p+\frac{(X_i-q_p)_+}{1-p}.
\end{equation}

\smallskip
\noindent\textbf{Expectiles.}
Fix $\tau\in(0,1)$ and put
\[
g_\tau(x,z):=\tau(x-z)_+-(1-\tau)(z-x)_+.
\]
The expectile $\xi_\tau$ is the unique solution of
$\mathbb E g_\tau(X_0,\xi_\tau)=0$. Its empirical counterpart
$\widehat\xi_{\tau,a:b}$ solves
$(b-a+1)^{-1}\sum_{i=a}^bg_\tau(X_i,\widehat\xi_{\tau,a:b})=0$.
If $F$ is continuous at $\xi_\tau$, define
\begin{equation}\label{eq:expectile-influence}
d_\tau:=\tau\{1-F(\xi_\tau)\}+(1-\tau)F(\xi_\tau),
\qquad
Y_i^\xi:=\frac{g_\tau(X_i,\xi_\tau)}{d_\tau}.
\end{equation}
Since $d_\tau\ge\min(\tau,1-\tau)>0$, the denominator is nonzero. At $\tau=1/2$, the expectile and its estimator reduce to the mean and sample mean. The influence function is well-known in expectile linearization, and we refer to \citet{holzmann:klar:2016} for the discussion of the i.i.d. setting.

\smallskip
\noindent\textbf{Lag-one linear prediction coefficient.}
Let $\mu=\mathbb E X_0$, suppose $\operatorname{Var}(X_0)>0$, and define
\begin{equation}\label{eq:ar-target}
\beta:=\frac{\mathbb E[(X_0-\mu)(X_1-\mu)]}{\operatorname{Var}(X_0)}.
\end{equation}
Then $\mu+\beta(X_{i-1}-\mu)$ is the best affine predictor of $X_i$ based on $X_{i-1}$ in mean squared error. A block estimator is the usual lag-one sample autocorrelation
\begin{equation}\label{eq:ar-estimator}
\widehat\beta_{a:b}
:=
\frac{\sum_{i=a+1}^b(X_{i-1}-\bar X_{a:b})(X_i-\bar X_{a:b})}
{\sum_{i=a}^b(X_i-\bar X_{a:b})^2},
\end{equation}
which is set to $0$ if we consider single-observation blocks (or when the empirical variance in the denominator vanishes). This estimator satisfies $|\widehat\beta_{a:b}|\le1$. An associated influence process is
\begin{equation}\label{eq:ar-influence}
Y_i^\beta
:=\frac{(X_{i-1}-\mu)(X_i-\mu)
-\beta(X_{i-1}-\mu)^2}{\operatorname{Var}(X_0)}.
\end{equation}
Notice that mathematically, the first influence term uses $X_{a-1}$, although the block estimator is calculated using only $X_a,\ldots,X_b$. This, however, is not a contradiction: in fact, this difference is simply absorbed in our theory by the block remainder. 
In an autoregressive model
\begin{equation}\label{eq:ar-model}
X_i-\mu=\phi(X_{i-1}-\mu)+u_i,\qquad |\phi|<1,
\end{equation}
we have $\beta=\phi$ whenever
$\mathbb E[(X_{i-1}-\mu)u_i]=0$. In this case we simply have 
$Y_i^\beta=(X_{i-1}-\mu)u_i/\operatorname{Var}(X_0)$.
The innovations $(u_i)$ may be dependent, provided this orthogonality condition holds. 

\textbf{Validation of Assumption \ref{ass:general}.} To demonstrate that the above examples really fit into our general theory, we need to impose low-level assumptions directly on the time series $(X_i)_i$. We give some sufficient conditions next.

\begin{assumption}
\label{ass:examples}
Let $r>4$. The observations $(X_i)_{i\in\mathbb Z}$ have a stationary causal representation as in Definition~\ref{def:causal}. For some $s>r$, $C<\infty$, and $\chi\in(0,1)$, assume
\begin{equation}\label{eq:example-moment-dependence}
\|X_0\|_s<\infty,
\qquad
\delta_s^X(k)\le C\chi^k,\qquad k\ge0.
\end{equation}
In addition, impose the following conditions for the specific parameter under consideration.

\begin{itemize}
\item[(i)] \textnormal{Quantile and expected shortfall.}
The marginal distribution has a differentiable density satisfying
\begin{equation}\label{eq:example-density}
\sup_{x\in\mathbb R}f(x)<\infty,
\qquad
\sup_{x\in\mathbb R}|f'(x)|<\infty,
\qquad f(q_p)>0.
\end{equation}

\item[(ii)] \textnormal{Expectile.}
We have $s>2r$, and $F$ has a bounded density.

\item[(iii)] \textnormal{Lag-one coefficient.}
We have $s>4r$ and $\operatorname{Var}(X_0)>0$.
\end{itemize}
Furthermore, we assume that the long-run variance $\sigma^2$ as defined in \eqref{e:LRV} is positive, for any  parameter to which the inference procedure is applied.
\end{assumption}

We can now validate our examples; the proof is given in the supplementary appendix.

\begin{theorem}[Verification of the high-level assumptions]
\label{thm:examples}
Suppose the respective conditions of Assumption~\ref{ass:examples} hold. Then, each of the (block) estimators is almost surely finite and admits the linearization \eqref{eq:method-linearization} with the corresponding centered stationary causal influence process. That process has a finite $r$th moment and geometrically decaying $L^r$ physical dependence coefficients. The block remainders satisfy
\begin{equation}\label{eq:example-remainder-rates}
\sup_{a\in\mathbb Z}\|R_{a:a+\ell-1}\|_r
\le
\begin{cases}
C\ell^{-3/4}\{\log(\mathrm e+\ell)\}^{3/2}, & \text{quantile},\\
C\ell^{-1}\{\log(\mathrm e+\ell)\}^2, & \text{expected shortfall},\\
C\ell^{-1}\log(\mathrm e+\ell), & \text{expectile},\\
C\ell^{-1}, & \text{lag-one coefficient},
\end{cases}
\qquad \ell\ge1.
\end{equation}
Consequently, Assumption~\ref{ass:general} holds with any
$\eta\in(1/r,1/4)$ for quantiles, and with any
$\eta\in(1/r,1/2)$ for the other parameters.
Theorem~\ref{thm:main} therefore applies to both normalizers.
\end{theorem}

The quantile bound follows from the probability estimates underlying \citet[Theorem~4]{wu:2005:quantiles}, together with moment control on the exceptional event. 
We also notice that in the lag-one example, our theory applies to \eqref{eq:ar-model} with non-Gaussian, physically dependent innovations, which covers a much larger class of data generating processes beyond standard AR(1). 

\section{Numerical experiments}
\label{sec:simulations}

We investigate the finite-sample performance of the proposed confidence
sequences and compare them with existing methods.

\subsection{Finite-sample coverage and interval widths}

We begin by assessing the simultaneous coverage of the two self-normalized
procedures and comparing the widths of the resulting confidence intervals.

\noindent\textbf{Data generation.}
We consider Gaussian and skewed data with dependence parameter
$\rho\in\{0,0.3,0.6\}$. In the Gaussian case, we generate a stationary
AR(1) process
\[
X_i^{(1)}
=\rho X_{i-1}^{(1)}+\sqrt{1-\rho^2}\,\varepsilon_i,
\qquad \varepsilon_i\stackrel{\mathrm{iid}}{\sim}N(0,1).
\]
For the skewed case, let $(Y_{i,k})_i$, $k=1,\ldots,5$, be independent
stationary Gaussian AR(1) processes satisfying
\[
Y_{i,k}
=\sqrt{\rho}\,Y_{i-1,k}
 +\sqrt{1-\rho}\,\varepsilon_{i,k},
\qquad \varepsilon_{i,k}\stackrel{\mathrm{iid}}{\sim}N(0,1),
\]
and set
\[
X_i^{(2)}
=\frac{\sum_{k=1}^{5}Y_{i,k}^2-5}{\sqrt{10}}.
\]
All autoregressions are initialized in their stationary distributions.
Both data-generating processes have mean zero and variance one. Their marginal distributions are
standard Gaussian and standardized $\chi^2_5$, respectively.
Importantly, the marginal distributions are the same for all $\rho$, which only measures serial dependence and where
 $\rho=0$ implies independence.

\medskip
\noindent\textbf{Simulation settings.}
We consider quantiles and expectiles at levels $0.5$ (median, mean) and $0.95$. For each marginal distribution and dependence
level, we generate 2{,}000 Monte Carlo replications, using the same paths
for both methods. As burn-in lengths we consider $m\in\{250,500,1000\}$.
For the three-block normalizer, we set $\kappa=1/4$.
The targeted simultaneous coverage is 95\%.
We use the weights
\[
\rho_{\mathrm{fw}}(x)=\log(\mathrm e+\log x),
\qquad
\rho_{\mathrm{LIL}}(x)=2\sqrt{\log(\mathrm e+\log x)},
\qquad x\ge1,
\]
for the forward and three-block methods, respectively.
The method-specific Brownian critical values are approximated once using
20{,}000 Brownian replications and held fixed across the data simulations.

\noindent\textbf{Evaluation.}
Although our theory concerns an unbounded monitoring horizon, monitoring needs to stop at some point. We therefore consider $n=m, m+1,\ldots,50{,}000$ in the simulations.
A replication counts as covered only if every confidence interval
throughout this period contains the true parameter.
Table~\ref{tab:sn-coverage} reports the resulting empirical simultaneous
coverage probabilities. Coverages are overall quite reasonable. Nominal approximations are better for $0.5$ quantiles/expectiles than for $0.95$ (as would be expected). At $m=250$, coverage for the upper quantile and expectile can fall to around $90\%$, particularly under stronger dependence. Overall coverage properties of the three--block normalizer are slightly better. It is however, also the computationally more expensive normalizer.

\begin{figure}[h!]
\centering
\begin{subfigure}[t]{0.49\linewidth}
  \centering
  \includegraphics[width=\linewidth]{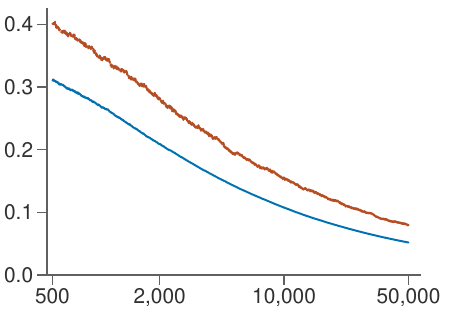}
  \caption{$m=500$}
  \label{fig:sn-width-500}
\end{subfigure}\hfill
\begin{subfigure}[t]{0.49\linewidth}
  \centering
  \includegraphics[width=\linewidth]{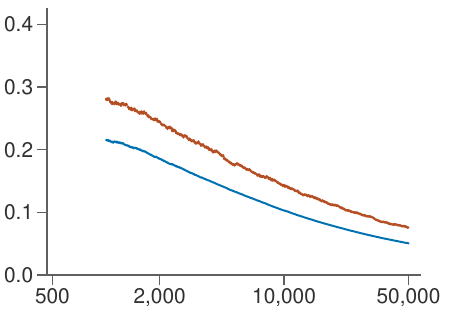}
  \caption{$m=1000$}
  \label{fig:sn-width-1000}
\end{subfigure}
\caption{Illustrative examples of pointwise median interval widths for $m=500, 1000$.  The parameter of interest is the  population median and we consider iid $N(0,1)$
data. Blue line: three-block; red line:
forward. The horizontal axis shows sample size $n$ on a logarithmic scale;
the vertical axis shows interval width. Total coverages are 95.30\% and 93.20\% in (a) and
94.15\% and 93.10\% in (b).}
\label{fig:sn-widths}
\end{figure}

We also compared interval widths, taking the median widths at each $n$ over all $2000$ Monte Carlo runs.
Figure~\ref{fig:sn-widths} provides two illustrative examples, showing
 widths for inference on the population median under iid
Gaussian sampling, with $m=500$ and $m=1000$.
In both examples, the curves exhibit a clear downward trend, with the
three-block method producing smaller median widths than the forward method
throughout the displayed monitoring period.

\begin{table}[h]
\centering
\caption{Simultaneous coverage (\%) of nominal 95\% confidence sequences over
$n=m,\ldots,50{,}000$, based on 2{,}000 paired replications.}
\label{tab:sn-coverage}
\begingroup
\small
\setlength{\tabcolsep}{3.2pt}
\renewcommand{\arraystretch}{0.96}
\begin{tabular*}{\linewidth}{@{\extracolsep{\fill}}llrrrrrr@{}}
\toprule
& & \multicolumn{2}{c}{$m=250$}
  & \multicolumn{2}{c}{$m=500$}
  & \multicolumn{2}{c}{$m=1000$} \\
\cmidrule(lr){3-4}\cmidrule(lr){5-6}\cmidrule(l){7-8}
$\rho$ & Parameter
& 3-block & Forward & 3-block & Forward & 3-block & Forward \\
\midrule
\multicolumn{8}{@{}l}{\itshape Gaussian} \\
\addlinespace[0.35ex]
$0.0$ & Mean & 94.55 & 94.05 & 94.75 & 95.20 & 94.00 & 94.00 \\
 & Median & 95.15 & 91.85 & 95.30 & 93.20 & 94.15 & 93.10 \\
 & $0.95$ quantile & 94.80 & 91.40 & 95.35 & 92.00 & 95.25 & 93.00 \\
 & $0.95$ expectile & 93.75 & 94.40 & 94.65 & 95.00 & 93.95 & 95.00 \\
\addlinespace[0.45ex]
$0.3$ & Mean & 94.90 & 94.50 & 95.30 & 94.50 & 95.00 & 94.65 \\
 & Median & 96.10 & 93.60 & 95.25 & 93.50 & 94.60 & 94.55 \\
 & $0.95$ quantile & 93.10 & 89.30 & 94.40 & 91.45 & 95.10 & 92.75 \\
 & $0.95$ expectile & 93.35 & 93.50 & 95.05 & 94.75 & 94.60 & 93.95 \\
\addlinespace[0.45ex]
$0.6$ & Mean & 94.35 & 94.10 & 94.75 & 94.25 & 95.25 & 95.05 \\
 & Median & 94.00 & 93.05 & 94.25 & 93.60 & 94.65 & 93.65 \\
 & $0.95$ quantile & 92.10 & 88.45 & 93.30 & 91.15 & 93.90 & 92.15 \\
 & $0.95$ expectile & 92.65 & 92.25 & 93.70 & 93.00 & 94.60 & 94.15 \\
\addlinespace[0.7ex]
\multicolumn{8}{@{}l}{\itshape Standardized $\chi^2_5$} \\
\addlinespace[0.35ex]
$0.0$ & Mean & 95.25 & 94.95 & 94.45 & 94.55 & 94.65 & 94.65 \\
 & Median & 95.45 & 92.25 & 94.65 & 93.30 & 95.35 & 93.55 \\
 & $0.95$ quantile & 94.40 & 90.15 & 95.00 & 91.65 & 95.05 & 92.55 \\
 & $0.95$ expectile & 93.50 & 93.60 & 94.45 & 94.85 & 94.60 & 94.30 \\
\addlinespace[0.45ex]
$0.3$ & Mean & 94.45 & 93.70 & 94.30 & 94.30 & 94.70 & 93.90 \\
 & Median & 94.90 & 92.55 & 94.70 & 93.55 & 94.60 & 93.05 \\
 & $0.95$ quantile & 92.75 & 89.85 & 93.85 & 91.45 & 95.60 & 93.15 \\
 & $0.95$ expectile & 91.80 & 91.70 & 94.20 & 93.50 & 94.75 & 94.45 \\
\addlinespace[0.45ex]
$0.6$ & Mean & 93.80 & 94.50 & 94.05 & 94.40 & 94.75 & 94.40 \\
 & Median & 94.50 & 93.25 & 95.55 & 94.05 & 94.55 & 94.10 \\
 & $0.95$ quantile & 91.45 & 89.20 & 92.60 & 90.00 & 93.75 & 91.95 \\
 & $0.95$ expectile & 89.65 & 90.85 & 91.70 & 91.45 & 92.85 & 93.50 \\
\bottomrule
\end{tabular*}
\par\vspace{3pt}
\begin{minipage}{\linewidth}
\footnotesize
$\rho$ denotes lag-one correlation. Three-block SN uses $\ell=\lfloor \kappa m\rfloor$ with $\kappa=1/4$.
The Monte Carlo standard error near 95\% is approximately 0.49 percentage points,
conditional on the calibrated critical values.
\end{minipage}
\endgroup
\end{table}


\subsection{Comparison with existing methods}
\label{sec:comparison-methods}

We next compare the proposed procedures with existing methods for
quantiles and the mean. Although our main interest lies in nonlinear
parameters, including the mean allows comparisons with additional
established benchmarks. For expectiles, expected shortfall and lag-one dependence parameters,
we did not find directly comparable benchmarks that fit our scenario.

Two distinctions are relevant to these comparisons. First, some methods assume
independence, whereas others are compatible with temporal dependence. Second,
confidence sequences provide simultaneous coverage throughout monitoring,
while fixed-width procedures guarantee coverage (only) at their prescribed
stopping time. Confidence sequences can also be stopped at a target width,
but typically require more observations to reach the same precision.

For these comparisons, we  present only results for Gaussian data;  corresponding results for skewed
data are in the supplementary appendix.

\medskip
\noindent\textbf{Quantiles.}
We consider two benchmarks.
\cite{howard:ramdas:2022} provide finite-sample confidence sequences for
population quantiles under iid sampling.
\cite{flegal:gong:2015} study width-based stopping procedures in the
context of MCMC, including estimation of quantiles. Their framework
allows dependent observations under regularity conditions and attains
nominal coverage at the prescribed stopping time as the precision
tolerance tends to zero.

\begin{figure}[!t]
\centering
\begin{subfigure}[t]{.49\linewidth}\centering
\includegraphics[width=\linewidth]{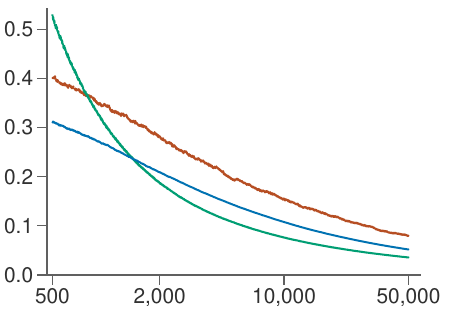}
\caption{$p=0.5$}\label{fig:comparison-q50}
\end{subfigure}\hfill
\begin{subfigure}[t]{.49\linewidth}\centering
\includegraphics[width=\linewidth]{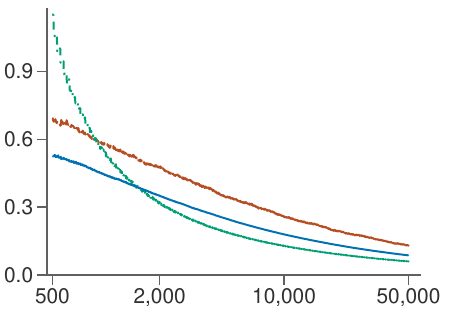}
\caption{$p=0.95$}\label{fig:comparison-q95}
\end{subfigure}
\caption{Illustrative pointwise median interval widths for quantiles of iid $N(0,1)$. Blue line: three-block; red line: forward; green line: Howard--Ramdas. The horizontal axis is sample size  and the vertical axis is interval width. Total coverages are 95.30/93.20/96.20\% in (a) and 95.35/92.00/96.30\% in (b).}
\label{fig:quantile-comparison-widths}
\end{figure}

Table~\ref{tab:quantile-stopping} illustrates the difference between the
two coverage guarantees. Flegal--Gong has the shortest stopping times,
with coverage at stopping close to $95\%$ for the median but around
$91\%$ for the upper quantile. Its simultaneous coverage is around
$60\%$, as expected for a procedure designed for a single prescribed
stop. Howard--Ramdas provides slightly conservative simultaneous
coverage under independence and reaches the target with about two-thirds
of the observations required by the self-normalized methods. Under
dependence, where its guarantee does not apply, coverage deteriorates
substantially and falls below $50\%$ for the median.
The self-normalized methods give much more stable coverage under
dependence. Three-block has consistently higher coverage than forward
at comparable stopping times; the difference is more pronounced for the
upper quantile, where forward coverage falls to about $91\%$.
Coverage at stopping is close to $100\%$ for both self-normalized methods,
and for Howard--Ramdas under independence. This reflects the price of
simultaneous guarantees: the wider intervals are conservative when used
for a single precision-based decision.
Figure~\ref{fig:quantile-comparison-widths} illustrates the corresponding
width differences. 

\medskip
\noindent\textbf{Mean.}
For the mean, we consider two additional benchmarks.
\cite{gnettner:kirch:2025} provide asymptotic confidence sequences that
accommodate weak dependence using a strongly consistent LRV estimator.
For \cite{waudby-smith:arbour:sinha:kennedy:ramdas:2024}, we use their iid
construction based on the sample variance. We additionally include
Flegal--Gong. For Gnettner--Kirch and Flegal--Gong, we use the same
nonoverlapping batch-means estimator with batch length
$\lfloor\sqrt{n}\rfloor$ under both independence and dependence.

\begin{figure}[!t]
\centering
\begin{subfigure}[t]{0.485\linewidth}
\centering
\includegraphics[width=\linewidth]{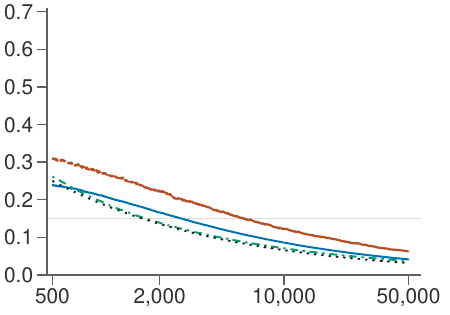}
\caption{Independent observations.}
\label{fig:mean-width-iid}
\end{subfigure}\hfill
\begin{subfigure}[t]{0.485\linewidth}
\centering
\includegraphics[width=\linewidth]{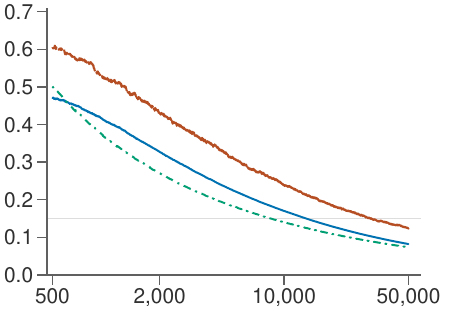}
\caption{Dependent observations, $\rho=0.6$.}
\label{fig:mean-width-dep}
\end{subfigure}
\caption{Illustrative pointwise median interval widths for the mean
under Gaussian observations. Blue line: three-block; red line: forward;
green line: Gnettner--Kirch; black line: Waudby-Smith et al.
The horizontal axis is sample size and the vertical axis is interval width.
Total coverages for three-block, forward, and Gnettner--Kirch are
94.75/95.20/93.10\% in (a) and 94.75/94.25/93.25\% in (b).
The method of Waudby-Smith et al.\ is shown only under independence, with total coverage
95.05\%.}
\label{fig:mean-widths}
\end{figure}

Table~\ref{tab:mean-comparison} and Figure~\ref{fig:mean-widths} show that
both self-normalized methods maintain coverage near the nominal level
under dependence. Waudby-Smith et al.'s method works well under independence but
substantially undercovers under dependence. Gnettner--Kirch reaches the
target considerably earlier than both self-normalized methods, although
its simultaneous coverage is around $93\%$.
The mean provides a favorable setting for LRV estimation: in these simple
models, and with $m=500$, standard batch means can estimate the LRV
directly from the observations. For nonlinear parameters this often
requires additional estimation steps. The comparison therefore illustrates
the cost of self-normalization when LRV estimation is relatively simple.

\begin{table}[!t]
\centering
\caption{Target-width stopping and coverage, based on 2{,}000 paired replications with $m=500$.}
\label{tab:quantile-stopping}
\begingroup\small
\setlength{\tabcolsep}{4pt}
\renewcommand{\arraystretch}{.96}
\begin{tabular*}{\linewidth}{@{\extracolsep{\fill}}llrrrr@{}}
\toprule
$\rho$ & Method & Median $T$ & Reached (\%)
& \shortstack{Simultaneous\\coverage (\%)}
& \shortstack{Coverage at\\stopping (\%)} \\
\midrule
\multicolumn{6}{@{}l}{\itshape Gaussian: $p=0.5$, $\epsilon=0.15$} \\
$0$ & Three-block    & 4{,}446  & 100.00 & 95.30 & 99.40 \\
    & Forward        & 4{,}268  &  99.95 & 93.20 & 98.50 \\
    & Howard--Ramdas & 2{,}794  & 100.00 & 96.20 & 99.95 \\
    & Flegal--Gong   & 1{,}500  & 100.00 & 60.60 & 93.55 \\
$0.6$ & Three-block    & 15{,}683 & 100.00 & 94.25 & 100.00 \\
      & Forward        & 17{,}888 &  93.20 & 93.60 & 99.70 \\
      & Howard--Ramdas & 2{,}766  & 100.00 & 46.35 & 93.55 \\
      & Flegal--Gong   & 3{,}500  & 100.00 & 60.85 & 93.30 \\
\addlinespace[.5ex]
\multicolumn{6}{@{}l}{\itshape Gaussian: $p=0.95$, $\epsilon=0.25$} \\
$0$ & Three-block    & 4{,}480  & 100.00 & 95.35 & 99.30 \\
    & Forward        & 4{,}297  & 100.00 & 92.00 & 97.10 \\
    & Howard--Ramdas & 2{,}811  & 100.00 & 96.30 & 99.90 \\
    & Flegal--Gong   & 1{,}000  & 100.00 & 57.15 & 91.85 \\
$0.6$ & Three-block    & 10{,}697 & 100.00 & 93.30 & 99.80 \\
      & Forward        & 11{,}560 &  97.00 & 91.15 & 98.05 \\
      & Howard--Ramdas & 2{,}806  & 100.00 & 71.30 & 97.25 \\
      & Flegal--Gong   & 2{,}000  & 100.00 & 56.80 & 90.70 \\
\bottomrule
\end{tabular*}
\par\vspace{3pt}
\begin{minipage}{\linewidth}\footnotesize
$\epsilon$ is the target full interval width; Reached is its attainment
rate by 50{,}000. Simultaneous coverage: all issued intervals through
50{,}000 contain the true quantile. Coverage at stopping: the single
interval at the stop, or at the cap if the target is not
reached ($T=\infty$ for the median), contains the true quantile. SN and Howard--Ramdas issue
intervals at every integer $n\ge500$; Flegal--Gong starts at 1000 and
checks every 500 observations.
\end{minipage}
\endgroup
\end{table}

\begin{table}[!t]
\centering
\caption{Mean inference: 2{,}000 paired replications per setting,
nominal coverage $95\%$, $m=500$, and target full width
$\epsilon=0.15$.}
\label{tab:mean-comparison}
\begingroup
\small
\setlength{\tabcolsep}{4pt}
\renewcommand{\arraystretch}{1.02}
\begin{tabular*}{\linewidth}{@{\extracolsep{\fill}}lrrrr@{}}
\toprule
Method & Median $T$
& \shortstack{Target reached\\(\%)}
& \shortstack{Simultaneous\\coverage (\%)}
& \shortstack{Coverage at\\stopping (\%)} \\
\midrule
\multicolumn{5}{@{}l}{\itshape Gaussian, $\rho=0$} \\
Three-block         & 2{,}284  & 100.00 & 94.75 & 98.00 \\
Forward             & 2{,}215  & 100.00 & 95.20 & 98.30 \\
Gnettner--Kirch     & 1{,}443  & 100.00 & 93.10 & 98.85 \\
Waudby-Smith et al. & 1{,}594  & 100.00 & 95.05 & 99.70 \\
Flegal--Gong        & 1{,}000  & 100.00 & 60.70 & 93.30 \\
\addlinespace[.4ex]
\multicolumn{5}{@{}l}{\itshape Gaussian, $\rho=0.6$} \\
Three-block         & 13{,}177 & 100.00 & 94.75 & 100.00 \\
Forward             & 14{,}689 &  95.70 & 94.25 & 99.90 \\
Gnettner--Kirch     & 7{,}396  & 100.00 & 93.25 & 99.90 \\
Waudby-Smith et al. & 1{,}588  & 100.00 & 23.50 & 87.00 \\
Flegal--Gong        & 3{,}000  & 100.00 & 61.70 & 93.00 \\
\bottomrule
\end{tabular*}
\par\vspace{3pt}
\begin{minipage}{\linewidth}
\footnotesize
$\epsilon$ is the target full interval width; Reached is its attainment
rate by 50{,}000. Simultaneous coverage: all issued intervals through
50{,}000 contain the true mean. Coverage at stopping: the single
interval at the stop, or at the cap if the target is not
reached ($T=\infty$ for the median), contains the true mean. Confidence sequences issue intervals
at every integer $n\ge500$; Flegal--Gong starts at 1000 and checks every
500 observations.
\end{minipage}
\endgroup
\end{table}

\FloatBarrier

\section{Discussion}
\label{sec:discussion}

We developed the first self-normalized confidence sequences for parameters of
weakly dependent time series. The construction avoids direct estimation
of the long-run variance. This is particularly helpful for nonlinear
statistics, where conventional variance estimation can require several
additional estimation steps with parameter tuning. Self-normalization instead uses fluctuations
of the parameter estimates themselves. The exact same principle applies across
a broad range of parameters.

The two proposed normalizers offer complementary advantages: Forward normalization
uses only the sequential parameter estimates and requires little
additional computation. It is therefore attractive when data access is
restricted or additional block estimates are expensive. Three-block
normalization offers greater stability over long horizons and accommodates
weights with the optimal LIL growth rate. Forward normalization requires
weights that grow (slightly) more quickly. The choice thus depends on both
computational considerations and compatibility with the desired weight
function.

The numerical results illustrate the usefulness of the approach under
serial dependence. For the mean, established methods benefit from
relatively straightforward long-run variance estimation. However, our main interest lies in nonlinear
parameters under dependence, where comparable sequential methods are
often scarce (e.g. for quantiles, expectiles, expected shortfall). Self-normalization provides effective inference in these
settings through a broadly applicable construction. 

A promising direction for future empirical research is inference for parameters
estimated through computationally expensive numerical procedures.
In such settings, forward self-normalization could provide sequential
uncertainty quantification directly from the evolving parameter estimates,
without requiring additional evaluations on subsamples.

\section*{AI use declaration}
We have made use of large language models (LLMs) in the course of this
project. We have used LLMs to draft and check parts of the proofs, to
write and debug simulation code, to help interpret finite-sample results,
and to edit the manuscript. We developed the statistical methodology
without using LLMs. We reviewed the AI-assisted material and take
responsibility for the final manuscript and code.
\clearpage
\normalem
\bibliographystyle{plainnat}
\bibliography{anytime_sn_references}

@article{alexopoulos:goldsman:mokashi:tien:wilson:2019,
  author  = {Alexopoulos, C. and Goldsman, D. and Mokashi, A. C. and Tien, K.-W. and Wilson, J. R.},
  title   = {{Sequest}: A sequential procedure for estimating quantiles in steady-state simulations},
  journal = {Operations Research},
  year    = {2019},
  volume  = {67},
  number  = {4},
  pages   = {1162--1183},
  doi     = {10.1287/opre.2018.1829}
}

@article{andrews:1991,
  author  = {Andrews, D. W. K.},
  title   = {Heteroskedasticity and autocorrelation consistent covariance matrix estimation},
  journal = {Econometrica},
  year    = {1991},
  volume  = {59},
  number  = {3},
  pages   = {817--858},
  doi     = {10.2307/2938229}
}

@misc{bibaut:kallus:lindon:2024,
  author        = {Bibaut, A. and Kallus, N. and Lindon, M.},
  title         = {Near-optimal non-parametric sequential tests and confidence sequences with possibly dependent observations},
  howpublished  = {arXiv preprint},
  year          = {2024},
  archivePrefix = {arXiv},
  eprint        = {2212.14411v5},
  primaryClass  = {stat.ME},
  url           = {https://arxiv.org/abs/2212.14411v5}
}

@article{darling:robbins:1967,
  author  = {Darling, D. A. and Robbins, H.},
  title   = {Confidence sequences for mean, variance, and median},
  journal = {Proceedings of the National Academy of Sciences of the United States of America},
  year    = {1967},
  volume  = {58},
  number  = {1},
  pages   = {66--68},
  doi     = {10.1073/pnas.58.1.66}
}

@misc{dette:graw:2025,
  author        = {Dette, H. and Graw, C.},
  title         = {{Gaussian} differential private bootstrap by subsampling},
  howpublished  = {arXiv preprint},
  year          = {2025},
  archivePrefix = {arXiv},
  eprint        = {2505.01197},
  primaryClass  = {stat.ML},
  url           = {https://arxiv.org/abs/2505.01197}
}

@article{dwork:roth:2014,
  author  = {Dwork, C. and Roth, A.},
  title   = {The algorithmic foundations of differential privacy},
  journal = {Foundations and Trends in Theoretical Computer Science},
  year    = {2014},
  volume  = {9},
  number  = {3--4},
  pages   = {211--407},
  doi     = {10.1561/0400000042}
}

@article{flegal:gong:2015,
  author  = {Flegal, J. M. and Gong, L.},
  title   = {Relative fixed-width stopping rules for {Markov} chain {Monte Carlo} simulations},
  journal = {Statistica Sinica},
  year    = {2015},
  volume  = {25},
  number  = {2},
  pages   = {655--675},
  doi     = {10.5705/ss.2013.209}
}

@misc{ge:liu:2026,
  author        = {Ge, M. and Liu, Y.},
  title         = {Time-uniform dependent central limit theory},
  howpublished  = {arXiv preprint},
  year          = {2026},
  archivePrefix = {arXiv},
  eprint        = {2610.29288},
  primaryClass  = {math.ST},
  url           = {https://arxiv.org/abs/2610.29288}
}

@article{gnettner:kirch:2025,
  author  = {Gnettner, F. and Kirch, C.},
  title   = {A new and flexible class of sharp asymptotic time-uniform confidence sequences},
  journal = {Statistics \& Probability Letters},
  year    = {2025},
  volume  = {226},
  pages   = {110462},
  doi     = {10.1016/j.spl.2025.110462}
}

@article{holzmann:klar:2016,
  author  = {Holzmann, H. and Klar, B.},
  title   = {Expectile asymptotics},
  journal = {Electronic Journal of Statistics},
  year    = {2016},
  volume  = {10},
  number  = {2},
  pages   = {2355--2371},
  doi     = {10.1214/16-EJS1173}
}

@article{howard:ramdas:2022,
  author  = {Howard, S. R. and Ramdas, A.},
  title   = {Sequential estimation of quantiles with applications to {A/B} testing and best-arm identification},
  journal = {Bernoulli},
  year    = {2022},
  volume  = {28},
  number  = {3},
  pages   = {1704--1728},
  doi     = {10.3150/21-BEJ1388}
}

@article{howard:ramdas:mcauliffe:sekhon:2021,
  author  = {Howard, S. R. and Ramdas, A. and McAuliffe, J. and Sekhon, J.},
  title   = {Time-uniform, nonparametric, nonasymptotic confidence sequences},
  journal = {The Annals of Statistics},
  year    = {2021},
  volume  = {49},
  number  = {2},
  pages   = {1055--1080},
  doi     = {10.1214/20-AOS1991}
}

@article{kiefer:vogelsang:2002,
  author  = {Kiefer, N. M. and Vogelsang, T. J.},
  title   = {Heteroskedasticity--autocorrelation robust standard errors using the {Bartlett} kernel without truncation},
  journal = {Econometrica},
  year    = {2002},
  volume  = {70},
  number  = {5},
  pages   = {2093--2095},
  doi     = {10.1111/1468-0262.00366}
}

@article{kiefer:vogelsang:2005,
  author  = {Kiefer, N. M. and Vogelsang, T. J.},
  title   = {A new asymptotic theory for heteroskedasticity-autocorrelation robust tests},
  journal = {Econometric Theory},
  year    = {2005},
  volume  = {21},
  number  = {6},
  pages   = {1130--1164},
  doi     = {10.1017/S0266466605050565}
}

@article{kiefer:vogelsang:bunzel:2000,
  author  = {Kiefer, N. M. and Vogelsang, T. J. and Bunzel, H.},
  title   = {Simple robust testing of regression hypotheses},
  journal = {Econometrica},
  year    = {2000},
  volume  = {68},
  number  = {3},
  pages   = {695--714},
  doi     = {10.1111/1468-0262.00128}
}

@article{lobato:2001,
  author  = {Lobato, I. N.},
  title   = {Testing that a dependent process is uncorrelated},
  journal = {Journal of the American Statistical Association},
  year    = {2001},
  volume  = {96},
  number  = {455},
  pages   = {1066--1076},
  doi     = {10.1198/016214501753208726}
}

@article{mueller:2007,
  author  = {M{\"u}ller, U. K.},
  title   = {A theory of robust long-run variance estimation},
  journal = {Journal of Econometrics},
  year    = {2007},
  volume  = {141},
  number  = {2},
  pages   = {1331--1352},
  doi     = {10.1016/j.jeconom.2007.01.019}
}

@article{shao:2010,
  author  = {Shao, X.},
  title   = {A self-normalized approach to confidence interval construction in time series},
  journal = {Journal of the Royal Statistical Society: Series B (Statistical Methodology)},
  year    = {2010},
  volume  = {72},
  number  = {3},
  pages   = {343--366},
  doi     = {10.1111/j.1467-9868.2009.00737.x}
}

@article{shao:2015,
  author  = {Shao, X.},
  title   = {Self-normalization for time series: A review of recent developments},
  journal = {Journal of the American Statistical Association},
  year    = {2015},
  volume  = {110},
  number  = {512},
  pages   = {1797--1817},
  doi     = {10.1080/01621459.2015.1050493}
}

@article{slavkovic:seeman:2023,
  author  = {Slavkovi{\'c}, A. and Seeman, J.},
  title   = {Statistical data privacy: A song of privacy and utility},
  journal = {Annual Review of Statistics and Its Application},
  year    = {2023},
  volume  = {10},
  pages   = {189--218},
  doi     = {10.1146/annurev-statistics-033121-112921}
}

@article{waudby-smith:arbour:sinha:kennedy:ramdas:2024,
  author  = {Waudby-Smith, I. and Arbour, D. and Sinha, R. and Kennedy, E. H. and Ramdas, A.},
  title   = {Time-uniform central limit theory and asymptotic confidence sequences},
  journal = {The Annals of Statistics},
  year    = {2024},
  volume  = {52},
  number  = {6},
  pages   = {2613--2640},
  doi     = {10.1214/24-AOS2408}
}

@article{wu:2005:dependence,
  author  = {Wu, W. B.},
  title   = {Nonlinear system theory: Another look at dependence},
  journal = {Proceedings of the National Academy of Sciences of the United States of America},
  year    = {2005},
  volume  = {102},
  number  = {40},
  pages   = {14150--14154},
  doi     = {10.1073/pnas.0506715102}
}

@article{wu:2005:quantiles,
  author  = {Wu, W. B.},
  title   = {On the {Bahadur} representation of sample quantiles for dependent sequences},
  journal = {The Annals of Statistics},
  year    = {2005},
  volume  = {33},
  number  = {4},
  pages   = {1934--1963},
  doi     = {10.1214/009053605000000291}
}

\end{document}